\documentclass[12pt, twoside,a4paper]{article}
\pdfoutput=1
\usepackage[cp1251]{inputenc}
\usepackage{ifthen}
\usepackage{amsmath}
\usepackage{amsfonts}
\usepackage{latexsym}
\usepackage{amsthm}
\usepackage{amssymb}
\newcommand{\CopyName}{ V.\ M.\ Zhuravlov}
\newcommand{\NAME}{ V.\ M.\ Zhuravlov}
\newcommand{\Year}{2024}
\newcommand{\rightheadtext}{Multiplicative logic in arithmetic}
     \newcounter{chapter}
     \newcounter{artpage}[chapter]
     \newcommand{\vs}{\vspace{.1in}}
     \newcommand{\vsk}{\vspace{.2in}}
\makeatletter
     \renewcommand{\@evenhead}{\footnotesize \ifthenelse{\value{artpage}=0}
     {\hfil}{\thepage\hfil \textsc {\leftmark} \hfil } }
     \renewcommand{\@oddhead}{\footnotesize\ifthenelse{\value{artpage}=0}
     {\hfil}{\hfil \textsc \rightmark \hfil \thepage} }
     \newcommand{\logo}{\baselineskip2pc \hbox to\hsize{\hfil\copyright\,\footnotesize
     \CopyName, \Year}}
     \renewcommand{\@oddfoot}{\ifthenelse{\value{artpage}=0}{\logo
     \refstepcounter{artpage}} {\hfil\refstepcounter{artpage}}}
     \renewcommand{\@evenfoot}{\ifthenelse{\value{artpage}=0}{\logo
     \refstepcounter{artpage}} {\hfil\refstepcounter{artpage}}}
     \renewcommand{\section}{\@startsection{section}{1}{0pt}{3.5ex plus
     1ex minus .2ex}{2.3ex plus 2.ex}{\large\hfil\textsc}}
     \vspace{3pt}\vs

\vs
\newcommand{\tit}{Multiplicative logic in arithmetic}
\date{2023}
\usepackage{hyperref}
\usepackage{graphicx}
\graphicspath{{Images/}}
\begin{document}
\hfill
\vspace{0.3in}
\markboth{{\NAME}}{{\rightheadtext}}\begin{center} \textsc {\CopyName} \end{center}\begin{center} \renewcommand{\baselinestretch}{1.3}\bf {\tit} \end{center}
\vspace{20pt plus 0.5pt} {\abstract{\noindent
The article explores the arithmetic of multiplication as a model of many-valued projective logic. It is demonstrated that closed numerical intervals within this framework constitute Heyting algebras. The conditions for these algebras to be Boolean are identified. The article's claims have undergone numerical verification. Paths for generalization to normed linear spaces are delineated.\newline
\textit{Multiplicative logic in arithmetic, 2024, msc: 03B50, 03B60;\vspace{3pt}}\newline
}
}\vsk
\tableofcontents
\section{Introduction}\par
This article serves as a continuation of my previous work:\href{https://arxiv.org/abs/2401.15117}{"Nonclassical logics and multivariate truth values."} In this piece, we explore a specific instance of projective logic (as defined in the aforementioned article), along with examples of logics that are not projective. These logics will be constructed based on the multiplicative properties of natural numbers. It is believed that such an approach could be beneficial for both logic and number theory, potentially even finding practical applications.\par
Let us first remind the reader that the \textbf{divisibility} relation (\textbf{\textit{x}} is a divisor of \textbf{\textit{y}}), being reflexive, antisymmetric, and transitive, constitutes an order relation on the set of natural numbers. This discrete order, \textbf{\textit{P}}, is strictly contained within the natural total order, \textbf{\textit{L}}, of the set of natural numbers, \textbf{\textit{N}}, induced by the successor operation: \textbf{\textit{(+1)}}; $ P \subset L \subset N \times N $. Unlike \textbf{\textit{L}}, divisibility does not form a linear order. By definition:\par
\textbf{Divisibility is a lattice order where the supremum of a pair of numbers is their least common multiple, and the infimum is their greatest common divisor.}\par
In the natural total order (where all numbers are comparable), the upper bound of a pair is their maximum, and the lower bound is their minimum. The mutual distributivity of these bounds is self-evident:
$$min(x,max(y,z)=max(min(x,y),min(x,z))$$
$$max(x,min(y,z)=min(max(x,y),max(x,z))$$
\textbf{A totally ordered set of natural numbers (as are all ordinals) is an unbounded above distributive discrete lattice with the least element} (it is understood that total order is equivalent to the existence of conjunctions for arbitrary sets of elements; and the scheme of axioms of natural induction determines that the truth values precisely form a set of natural numbers. However, these particularities are not of concern to us at this moment).
What are the logical consequences of this well-known fact?
\section{General and projective logic of natural truth values}\par
As per the previously mentioned article, \href{https://arxiv.org/pdf/2401.15117}{https://arxiv.org/pdf/2401.15117},—\textbf{ natural numbers are many-valued logic without "absolute truth" and without negation, but with a false value (that is 1) and with an infinite gradation of approximations to the truth.} Accordingly, the lattice dual to natural numbers will be logic that does not have a false meaning. This is a non-local part of natural number logic. And logics based on other linearly ordered sets will also have similar logical structures (basic operations — conjunction and disjunction, polysemy of truth values in the absence or incomplete presence of limiting values — truth and falsity). But there are also local properties of such logic — we talked about them in the article mentioned above.\par
According to our theorem on projective logic, \textbf{every closed natural interval} $Q_{a,b}\equiv\{x\in N|[a,b]\subset N\wedge(a\leq x\leq b)\}$ \textbf{is a Heyting algebra with zero \textit{a} and one \textit{b}} (this is a very special Heyting algebra — in it: $(x\longmapsto y)=b$ if $(x\leq y)$; and $(x\longmapsto y)=y$ if $(y\leq x)$; and therefore: $\neg x=a$ and $\neg\neg x=b$; this is how we prove the existence of Heyting pseudocomplements in an arbitrary lattice $Q_{a,b}$).\par
"What are the logical properties of natural numbers when ordered by divisibility? Let's recall another well-known fact. Denote by \textbf{NN} the set of all finite sequences of numbers from $\{0\}\bigcup N$. Then we have \textbf{an ordinal (and even a lattice) isomorphism:} $ N^* \longleftrightarrow NN $, \textbf{where} $N^{*}$ \textbf{are the natural numbers ordered by divisibility.} This isomorphism is established by decomposing any number \textbf{x} into powers of prime factors:
$$ x \longleftrightarrow \prod_{j=1}^{n} p_j^{x_j} $$
where $p_{j}$ are the first \textbf{n} prime numbers raised to some powers from an arbitrary finite sequence $x_{j}$ (zero powers are also possible; in other words, every number is a product of certain powers of prime numbers from the sequence of such an \textbf{n}-tuple of the first prime numbers, where all subsequent primes will only have a zero power). However, for our further reasoning, it would be better to standardize the multiplicative representation of numbers by eliminating the number \textbf{n}, which plays only a secondary role. As \textbf{NN}, we will consider infinite sequences from $\{0\}\bigcup N$ that contain only a finite number of non-zero elements:"
$$\forall j(j\in N,x_{j}\in\{0\}\bigcup N)\exists n(n\in N):\sum_{j\in N}(x_{j})\leq n$$
 — we could have written any other normalizing condition equivalent to the one stating that there are only a finite number of non-zero elements in the sequence: $(x_{j}=0)$ "for almost all" \textbf{j} (in fact, this concerns an ultrafilter of the finite subsets of the natural numbers). Then the aforementioned isomorphism can be described by the formula:
$$ x \longleftrightarrow \prod_{j \in \mathbb{N}} p_j^{x_j} $$
It is clear that \textbf{y} divides \textbf{x} completely (denoted by $y \prec x$) if and only if, in the above prime factor decomposition, it holds that $\forall j (y_j \leq x_j)$. Furthermore, this correspondence is a bijective epimorphism, as every number has a unique decomposition into powers of prime factors and vice versa. This proves the isomorphism between natural numbers and normalized infinite sequences of natural numbers. This information is sufficient to prove the following statements about the lower and upper bounds of the $\prec$-order (greatest common divisor and least common multiple).\newline
\textbf{The upper and lower bounds of the} $\prec$\textbf{-order possess the following properties:\newline
1.Idempotency\newline
2.Commutativity\newline
3.Associativity\newline
4.Mutual distributivity\newline
Which, obviously, once again demonstrates that divisibility is an order relation.}\newline
Proof. Every natural number has a canonical decomposition into prime factors, as indicated in the formula given above; which, as we have already mentioned, leads to an isomorphism between the lattice of numbers by their divisibility and the set of finite sequences of numbers, ordered in the "usual" way. Specifically, we represent any finite sequence as infinite by padding it with zeros. Then one number is a divisor of another if and only if the powers of the prime factors in its decomposition are less than or equal to the corresponding powers in the decomposition of the other number. This obviously proves points \textbf{1) - 3)}. Point \textbf{4)} follows from the previously mentioned formulas for the mutual distributivity of maxima and minima in the natural order of natural numbers. The proof is complete.\newline
\textbf{Note:} The first three points are well-known; as for the fourth, I have not encountered any mention of it; I believe this is because it was not considered important enough (just as the mutual distributivity of the maxima and minima of numbers was not considered important).
\section{Localization of truth, falsity, and complements in multiplicative logic}\par
As we have already demonstrated in the previously referenced article, within distributive multi-valued logic, any closed interval $Q_{y,x}$ constitutes a logic with absolute truth \textbf{x} and falsity \textbf{y} (given: $y\leq x$). Natural multiplicative logic is not an exception to this. However, what is the state of complements within this logic?\newline
\textbf{Theorem.} $Q_{y,x}$ \textbf{\textit{represents a finite (and thus discrete) Heyting algebra. The local complement of any given element \textbf{a} is the element}} $(\neg a)$\textbf{\textit{, which is calculated according to the formula:}}
$$\forall j:[(a_{j}>y_{j})\Longrightarrow((\neg a)_{j}=y_{j})]\wedge[(a_{j}=y_{j})\Longrightarrow((\neg a)_{j}=x_{j})]$$
\textbf{\textit{In general, the implication of}} $a\longrightarrow b$ \textbf{\textit{in this algebra will be the local complement of a in the algebra}} $Q_{(b\wedge a),x}$ \textbf{\textit{— i.e. the implication is determined by the formula:}}
$$\forall j:[(a_{j}>b_{j})\Longrightarrow((a\longrightarrow b)_{j}=b_{j})]\wedge[(a_{j}\leq b_{j})\Longrightarrow((a\longrightarrow b)_{j}=x_{j})]$$
\textbf{Note:} We refer to the usual Heyting pseudocomplement as the local complement to underscore its reliance on the endpoints of the closed interval.\newline\newline
\textbf{Proof.} In constructing the isomorphism $ N^* \longleftrightarrow NN $, we presented a series of known facts. Now, it only remains to apply these facts to the definitions from Heyting algebra. The negation $(\neg a)$ is the greatest element that is disjunctive with \textbf{a}. Therefore, in the factorization into prime factors, if the "coordinate" of \textbf{a} (that is, the power of the prime number) at a given prime number matches the coordinate of zero (i.e., \textbf{y}), then the coordinate of negation should be equated to the coordinate of one (i.e., \textbf{x}). Otherwise, the coordinate of a will be strictly greater than the coordinate of \textbf{y}, and then the coordinate of negation will be the maximal possible within the local algebra $ Q_{y,x} $, that is, "unitary" (remembering that $ x \geq a \geq y $).\par
The corresponding property of implication, a natural extension of negation, is proven similarly. Here, we discuss $ a \rightarrow b $ as the largest element whose intersection with \textbf{a} is less than \textbf{b}. As with negation, this fact is expressed in the coordinates of the decomposition into prime factors. \textbf{The proof is concluded.}\newline\newline
\textbf{Corollary 1. \textit{The algebra $ Q_{y,x} $ is Boolean if and only if the degree of each prime number in the factorization of y is less than or equal to the degree of this number in the factorization of x by no more than one.}}\newline
\textbf{Proof.} Let $ x \geq a \geq y $ and let $ (\neg a) \equiv (a \to y) $ be the complement to \textbf{a} in the algebra $ Q_{y,x} $; then, according to the corollary's condition, the degree of a prime number in the factorization of the number $(\neg a)$ will be equal to its degree in \textbf{x} if, and only if, the degree of this number in the factorization of \textbf{a} is the same as in \textbf{x}. Otherwise, this degree will be one less, that is, it will coincide with the degree of this prime number in the factorization of \textbf{y}. It follows that $ a \vee (\neg a) = x $; thus, the law of excluded middle is satisfied. It is evident that any element \textbf{a} from $ Q_{y,x} $ will be completely determined by the set of prime numbers which, in its factorization, have a degree equal to the degree of this prime number in the factorization of \textbf{x} (and then its complement will be determined by the corresponding degree in the factorization of \textbf{y}). It is clear that this argument is fully reversible. \textbf{The proof is concluded.}\par
\textbf{Note.} In the simplest case, all elements from $ Q_{y,x} $ will decompose only into the first (or zero) degrees of prime numbers. In the general case, we obtain a Boolean algebra that is isomorphic to the algebra of this simplest case, where all degrees are "shifted" by the same amount.\newline
\textbf{Corollary 2. The algebra} $ Q_{y,x} $ \textbf{is Boolean if and only if the local complement is calculated by the formula:}
$$ (\neg a) = \frac{x \cdot y}{a} $$
\textbf{Proof.} Indeed, according to the proof of Corollary 1,
$$ (a_j = x_j) \Rightarrow ((\neg a)_j = y_j) $$
$$ (a_j = y_j) \Rightarrow ((\neg a)_j = x_j) $$
— and this is equivalent to the formula of Corollary 2, after carrying out factorizations into powers of prime factors and simplifications. \textbf{The proof is concluded.}
\section{Computer verification of logical properties of multiplicative arithmetic.}\par
We will simply confirm these properties by examining a certain number of cases. Additionally, we will present a series of formulas and programs that further clarify the aforementioned processes. We will be using the \textbf{MATHCAD 15} software. Some of the formulas and programs are accompanied by specific examples.\par
Verification of the idempotency, commutativity, associativity, and mutual distributivity of the operations of the greatest common divisor (\textbf{gcd}) and the least common multiple (\textbf{lcm}):\begin{center}
\includegraphics{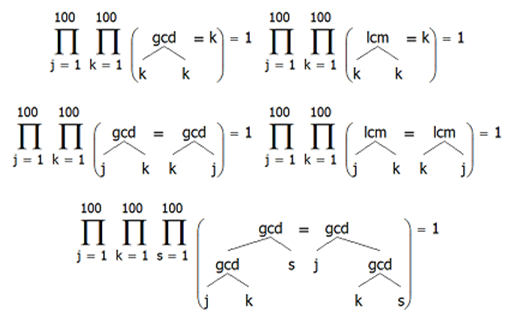}\end{center}\begin{center}
\includegraphics{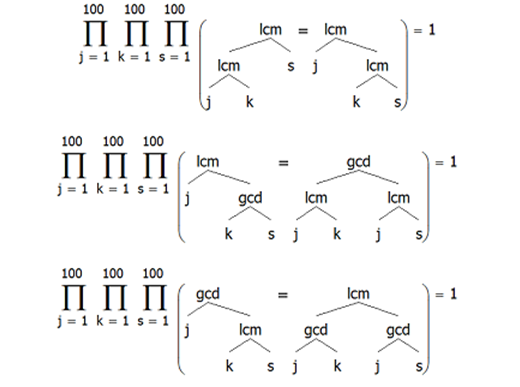}\end{center}
— all cases for numbers from 1 to 100 have been examined; there are no apparent fundamental reasons for these properties to fail at other number values.\par
Thus, the natural numbers form a distributive lattice with a minimum element (with respect to \textbf{gcd} — the greatest common divisor, and \textbf{lcm} — the least common multiple). Consequently, for each number, there are disjoint elements — any that are mutually prime with it. The divisibility relation establishes the order in this lattice. Furthermore, by taking any pair of numbers, we can construct a 4-element Boolean algebra, where \textbf{gcd} is the Boolean zero and \textbf{lcm} is the Boolean one (the algebra is 2-element if one member of the pair is a divisor of the other).\par
Here is the vector of the first \textbf{17} prime numbers:\begin{center}\includegraphics{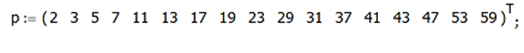}\end{center}
Translate the vector into the products of the corresponding powers of prime numbers, according to the above formula of isomorphism:\begin{center}\includegraphics{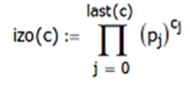}\end{center}
From vectors for powers of 3-numbers we make a complement vector:\begin{center}\includegraphics{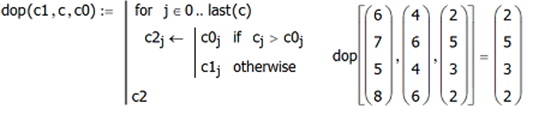}\end{center}\par
Finally. If $x\geq a\geq y$, then the negation of \textbf{a} is constructed as follows. Its simple expansion has powers equal to the powers of \textbf{y} — where the degrees of expansion of \textbf{a} are greater than the powers of \textbf{y}; where they are equal, we set these powers equal to the corresponding powers of \textbf{x}. This is implemented programmatically as follows (we first create a random input consisting of 3 integer 4-vectors, all coordinates of which do not increase from vector to vector; these will be powers of the first four prime numbers; the products of such powers will make up 3 numbers ordered by division; and then, we create a vector of degrees as close as possible to the “maximum” vector, but not exceeding the “minimum” one — where the average vector is greater than the minimum one and the addition in the Heyting algebra is constructed from it):\begin{center}\includegraphics{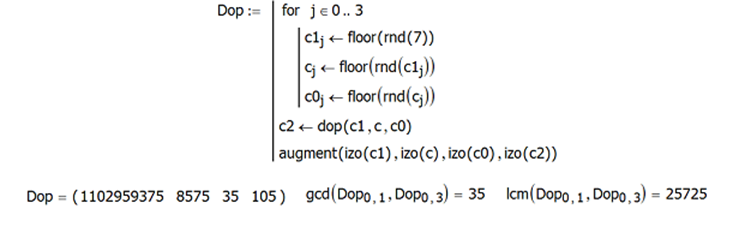}\end{center}\par
The Dop program specifies random (in a certain range) values of the powers of the first four prime numbers, satisfying the condition that the 2-nd number falls into the divisibility interval limited by the interval of the 1-st and 3-rd numbers. And finds the 4th number (from the same segment), locally complementary to the 2nd. A specific demonstration is shown. These demonstrations can be repeated indefinitely.\par
In the case when the degrees in the expansions of \textbf{x} and \textbf{y} differ from each other by no more than \textbf{1}, this is equivalent to the local complement formula: $(x\ast y)/a$ for $x\geq a\geq y$; and if this is so, then we have the Boolean algebra of the interval $Q_{y,x}$. Otherwise — Heyting's.
But we still need to define the implication in the Heyting algebra $Q_{y,x}$ and prove its properties. Namely, if: $x\geq a\geq y$ and $x\geq c\geq y$, then the implication between \textbf{a} and \textbf{c} in the algebra $Q_{y,x}$ is the complement to a in the algebra $Q_{y\wedge c,x}$. Those, — we replace a with c where a is greater than \textbf{c}, and put \textbf{x} in all other places; you need to write an algorithm and check it for implication properties.
Implication $c\longrightarrow c2$, when \textbf{c} and \textbf{c2} lie in the same interval:\begin{center}\includegraphics{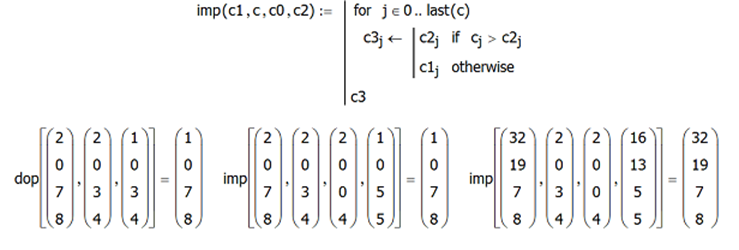}\end{center}
— the operation of the program is again illustrated with a specific numerical implementation.\newline
Further:\begin{center}\includegraphics{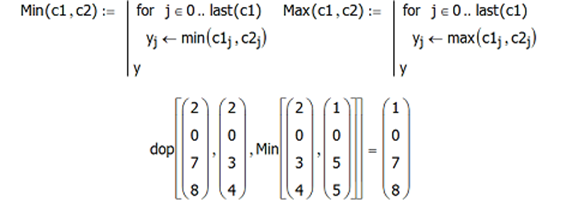}\end{center}
— from here we see that in the implication the role of zero is played by the infimum of both implied ones; Is it possible that this is the case for all Heyting algebras? In general, it turns out that the implication depends only on the 1-st we have chosen, but not on the 0-th.
It is known that the operation of implication $\longrightarrow$ is such that $(a\wedge b)\leq c$ is equivalent to $a\leq(b\longrightarrow c)$. Next we check this fact:\begin{center}\includegraphics{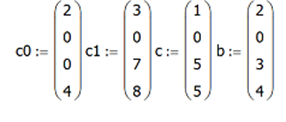}\end{center}\begin{center}\includegraphics{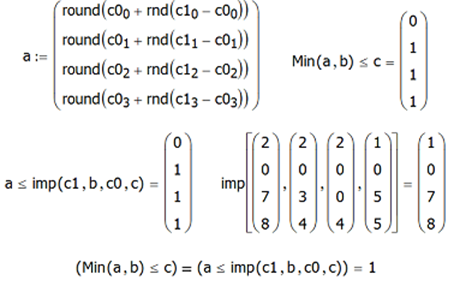}\end{center}
So, implication has one of the basic Heyting properties.
Let's also check one of the basic relations of projective logic: $x\wedge(z\vee y)=(x\wedge z)\vee y$ for arbitrary \textbf{z} and $x\geq y$:\begin{center}\includegraphics{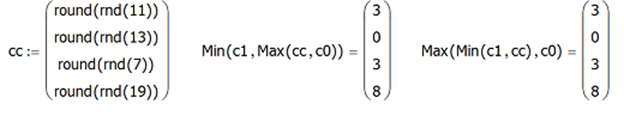}\end{center}
\section{Conclusion}\par
In the previous article (cited at the outset), we demonstrated that the most general model of many-valued logic should be considered an unbounded modular lattice with local complements (or equivalently, with local implications). If this lattice is distributive, then the local algebras of its closed intervals $Q_{y,x}$ will be Heyting algebras.\par
In this paper, we examined a particular case of such a lattice—the multiplicative set of natural numbers, ordered by the relation of divisibility. The upper bound of a pair of numbers is their least common multiple, and the lower bound—their greatest common divisor. Given that this lattice is isomorphic to the lattice of infinite sequences, almost everywhere zero, the following generalization suggests itself. If one takes a Hilbert (or Banach) space with a certain fixed (countable) basis, then its normalized vectors form a distributive unbounded lattice with respect to the operations of coordinate-wise minima and maxima of vector pairs. This lattice will also be a \textbf{local} Heyting algebra (continuous, as opposed to the natural divisibility lattice). It appears that if the basis is not fixed, we obtain a modular lattice. This once again brings us closer to the ideas of quantum logic.\par
Another path to generalization may be suggested by monoid theory. The divisibility relation in an arbitrary monoid is a pre-order relation. And if in the monoid there exist greatest common divisors and least common multiples for pairs of elements, we can construct a logic similar to the logic constructed in this article.\par
Moreover, all these (local with respect to the choice of zeros and ones) "interval" Heyting algebras are not just a set of algebras but are 'embedded' within a certain global logic (in which there is no ultimate truth or falsehood, but there is a multi-valued set of truth values). Thus, it is possible to construct general, global models of theories, along with their more specific local modifications.
\section{References}
\noindent [1]  Haskell B. Curry, McGraw-Hill Book Company, INC.\newline
\textit{Foundations of Manthematical logic}, (1969).\newline
[2] Garrett Birkhoff, Providence Rhode Island,\newline
\textit{Lattice theory}, (1967).\newline
[3] Helena Rasiowa and Roman Sikorski,  Panstwowe Wydawnlctwo Naukowe Warszawa,\newline
\textit{The Mathemayics of Metamathematics}, (1963).\newline
[4]Roman Sikorski, Springer-Verlag,\newline
\textit{Boolean Algebras}, (1964).\newline
[5] R. Goldblatt. North-Holland Publishing Company: Amsterdam, New York, Oxford.\newline
\textit{Topoi. The categorial analysis of logic}, (1979).\newline
\end{document}